\documentclass[smallextended,numbook,runningheads]{svjour3}     % onecolumn (second format)
\smartqed  % flush right qed marks, e.g. at end of proof
\usepackage{amsfonts,anysize,amsmath,indentfirst,geometry,xcolor}
\usepackage[colorlinks,citecolor=cyan,linkcolor=blue]{hyperref}
\usepackage{amssymb}
\usepackage{mathrsfs}
\usepackage[figuresleft]{rotating}
\usepackage{graphicx}
\usepackage{amsbsy}
\usepackage{eqnarray}
\usepackage{booktabs}
\usepackage{mathptmx}
\usepackage{multirow}
\usepackage{epstopdf}
\usepackage{array}
\usepackage{algorithm}
\usepackage{algorithmic}
\usepackage{graphicx}
\usepackage{wasysym}
\usepackage[misc]{ifsym}
\usepackage{cite}
\usepackage{epstopdf} % needed if you have eps figures in a pdflatex manuscript
\journalname{Numerical Algorithms}
\begin{document}

\title{Partial condition number for the equality constrained linear least squares problem\thanks{The work is supported by the National Natural Science Foundation of China (No. 11201507), the Fundamental Research Funds for the Central Universities (No. 106112015CDJXY100003), and the China Scholarship Council.}}


\author{Hanyu Li\and Shaoxin Wang  %etc.
}


\institute{H. Li (\Letter) \at College of Mathematics and Statistics, Chongqing University, Chongqing, 401331, P. R. China\\
           \email{lihy.hy@gmail.com or hyli@cqu.edu.cn }
           \and
                      S. Wang \at College of Mathematics and Statistics, Chongqing University, Chongqing, 401331, P. R. China \\
%              \at College of Mathematics and Statistics, Chongqing University, Chongqing, 401331, P. R. China \\
%              Tel.: +123-45-678910\\
%              Fax: +123-45-678910\\
              \email{shaoxin.w@gmail.com}           %  \\
}

\date{Received: date / Accepted: date}

\maketitle

\begin{abstract}
In this paper, the normwise condition number of a linear function of the equality constrained linear least squares solution called the partial condition number is considered. Its expression and closed formulae are first presented when the data space and the solution space are measured by the weighted Frobenius norm and the Euclidean norm, respectively. Then, we investigate the corresponding structured partial condition number when the problem is structured. To estimate these condition numbers with high reliability, the probabilistic spectral norm estimator and the small-sample statistical condition estimation method are applied and two algorithms are devised. The obtained results are illustrated by numerical examples.
%Include keywords and mathematical subject classification numbers as needed.
\keywords{Linear least squares problem \and Equality constraint \and Partial condition number \and Probabilistic spectral norm estimator \and Small-sample statistical condition estimation}
\subclass{65F20\and 65F35\and 65F30\and 15A12\and 15A60}
\end{abstract}

\section{Introduction and preliminaries}
The equality constrained linear least squares problem can be stated as follows:
\begin{equation}\label{1.1}
{\rm LSE:}\quad\mathop {\min }\limits_{Bx=d } \left\|b - Ax\right\|_2,
\end{equation}
where $A \in {\mathbb{R}^{m \times n}}$ and $B \in {\mathbb{R}^{s \times n}}$ with $m+s\geq n\geq s\geq0$, $b \in \mathbb{R}^{m}$ and $d \in \mathbb{R}^{s}$.
Hereafter, the symbols $\mathbb{R}^{m \times n}$ and $\mathbb{R}^{n}$ stand for the set of $m\times n$ real matrices and the real vector space of dimension $n$, respectively. To ensure that the LSE problem \eqref{1.1} has a unique solution, we need to assume that \cite{Bjo}
\begin{equation}\label{1.2}
{\rm rank}(B)=s,\quad {\rm null}(A)\cap{\rm null}(B)=\{0\}.
\end{equation}
The first condition in \eqref{1.2} implies that the linear system $Bx=d$ is consistent and hence that the LSE problem \eqref{1.1} has a solution, and vice versa; The second one, which says that the matrix $[A^{T}, B^{T}]^{T}$ is full column rank, guarantees that there is a unique solution to \eqref{1.1},  and vice versa. Here, for a matrix $C$, $C^{T}$ denotes its transpose. Throughout this paper, we assume that the conditions in \eqref{1.2} always hold.
In this case, the unique solution to the LSE problem \eqref{1.1} can be written as \cite{Bjo,Cox}
\begin{equation}\label{1.3}
x(A,B,b,d) = (AP)^ \dag  b +B_A^ \dag  d,
\end{equation}
where
\begin{equation*}\label{1.4}
P = I_n  - B^ \dag  B,\quad B_A^ \dag   = (I_n  - (AP)^ \dag  A)B^ \dag
\end{equation*}
with $I_n$ being the identity matrix of order $n$ and $ B^ \dag$ being the Moore-Penrose inverse of $B$. When $s=0$, i.e., $B=0$ and $d=0$, the LSE problem \eqref{1.1} reduces the classic linear least squares (LLS) problem
\begin{equation}\label{1.100}
{\rm LLS:}\quad\mathop {\min }\limits_{x\in\mathbb{R}^{n} } \left\|b - Ax\right\|_2,
\end{equation}
the conditions in \eqref{1.2} reduce to $ A$ being full column rank which ensures that the solution to \eqref{1.100} is unique, and the solution \eqref{1.3} reduces to $x(A,b) = A^ \dag  b$.

The LSE problem finds many important applications in some areas. For example, we will encounter it in the analysis of large scale structures, in signal processing, and in solving inequality constrained least squares problem \cite{Bar88b, Bjo, Law95}. So, some scholars considered its algorithms and perturbation analysis (see e.g., \cite{Bar88b, Bjo, Law95,Eld,Wei922,Cox}). An upper bound for the normwise condition number of the LSE problem was presented in \cite{Cox}, and the mixed and componentwise condition numbers and their easily computable upper bounds of this problem can be derived from \cite{Li14} as the special case.

In this paper, we mainly consider the partial condition number of the LSE problem when the data space $\mathbb{R}^{m \times n}\times \mathbb{R}^{s \times n}\times \mathbb{R}^{m}\times \mathbb{R}^{s}$ and the solution space $\mathbb{R}^{n}$ are measured by the weighted Frobenius norm
\begin{equation}\label{1.5}
\left\| (\alpha_A A, \alpha_B B,\alpha_b b,\alpha_d d) \right\|_F=\sqrt{\alpha_A^2\left\| A\right\|_F^2+\alpha_B^2\left\| B\right\|_F^2+\alpha_b^2\left\|b \right\|_2^2+\alpha_d^2\left\| d \right\|_2^2}
\end{equation}
with $ \alpha_A>0, \alpha_B>0, \alpha_b>0$, and $\alpha_d>0$, and the Euclidean norm $\left\|x \right\|_2$, respectively. As mentioned in Abstract, the partial condition number which is also called the subspace condition number \cite{Cao} is referred to the condition number of a linear function of the LSE solution $x(A,B,b,d)$, i.e., $L^{T}x(A,B,b,d)$ with $L\in\mathbb{R}^{n \times k}$ ($k\leq n$). This kind of condition number has some advantages. For example, when $L$ is the identity matrix or a column vector of the identity matrix, the partial condition number will reduce to the condition number of the solution $x(A,B,b,d)$ or of an element of the solution. Cao and Petzold first proposed the partial condition number for linear systems \cite{Cao}. Later, it was proposed for LLS problem and total least squares problem \cite{Ari, Bab}. In \cite{Cao,Ari, Bab}, the authors also provided some specific motivations for investigating this kind of condition number.

The idea on the weighted Frobenius norm can be traced back to Gratton \cite{Grat}, who derived the normwise condition number for the LLS problem \eqref{1.100}
based on the following weighted Frobenius norm
\begin{equation}\label{1.6}
\left\| (\alpha A, \beta b) \right\|_F=\sqrt{\alpha^2\left\| A\right\|_F^2+\beta^2\left\| b \right\|_2},\quad \alpha>0, \beta>0.
\end{equation}
Subsequently, this kind of norm was used for the partial condition number for the LLS problem \cite{Ari} and the normwise condition number of the truncated singular value solution of a linear ill-posed problem \cite{Ber}. As pointed out in \cite{Grat}, this norm is very flexible. With it, we can monitor the perturbations on $A$ and $b$. For example, if $\alpha\rightarrow \infty$, no perturbation on $A$ will be permitted; similarly, if $\beta\rightarrow \infty$, there will be no perturbation on $b$ allowed. Obviously, the norm in \eqref{1.5} is a simple generalization of the one in \eqref{1.6}, and is also very flexible. There is another kind of generalization of the norm in \eqref{1.6}: $\left\| (T A, \beta b) \right\|_F$, which was used  by Wei et al. in \cite{Wei} for the normwise condition number of rank deficient LLS problem. Here, $T$ is a positive diagonal matrix.

Like the structured linear systems and the structured LLS problem, the structured LSE problem arises in many applications, e.g., in signal processing and the area of optimization \cite{Bjo, Law95}. Rump \cite{Rump1,Rump2} presented the perturbation theory for the structured linear systems with respect to normwise distances and componentwise distances. The obtained results generalized the corresponding ones in \cite{Hig}. For the structured LLS problems, Xu et al. \cite{Xu} considered their structured normwise condition numbers, and Cucker and Diao \cite{Cucker} presented their structured mixed and componentwise condition numbers. In addition, the structured condition numbers for the total least squares problem were provided by Li and Jia in \cite{Lijia}. The results in \cite{Lijia,Rump1,Rump2,Xu} show that the structured condition number can be much tighter than the unstructured one. So, based on the study on the partial condition number, we also investigate the structured partial condition number of the structured LSE problem.

The rest of this paper is organized as follows. Section 2 presents the expression and closed formulae of the partial condition number for the LSE problem. The expression of the corresponding structured partial condition number is given in Section 3. On basis of the probabilistic spectral norm estimator by Hochstenbach \cite{Hochs13} and the small-sample statistical condition estimation (SSCE) method by Kenney and Laub \cite{Kenney94}, Section 4 is devoted to the statistical estimates and algorithms of the results derived in Sections 2 and 3. The numerical experiments for illustrating the obtained results are provided in Section 5. %Finally, we present the concluding remarks of the whole paper.

Before moving to the following sections, we first introduce some results on the operator 'vec' and Kronecker product, and the generalized singular value decomposition (SVD) of a matrix pair. They will be necessary later in this paper.

For a matrix $A=[a_1,\cdots, a_n]\in
\mathbb{R}^{m\times n}$ with $a_i\in \mathbb{R}^{m}$, the operator 'vec' is defined as follows
\begin{equation*}
{\rm vec}(A)=[a_1^{T},\cdots,a_n^{T}]^{T}\in \mathbb{R}^{mn}.
\end{equation*}
Let $A =(a_{ij})\in {\mathbb{R}^{m \times n}}$ and $B \in
{\mathbb{R}^{p \times q}}$. The {\em Kronecker product} between $A$ and $B$ is defined
by (see,
e.g., \cite[Chapter 4]{Horn91})
\[A \otimes B = \left[ {\begin{array}{*{20}{c}}
{{a_{11}}B}&{{a_{12}}B}& \cdots &{{a_{1n}}B}\\
{{a_{21}}B}&{{a_{22}}B}& \cdots &{{a_{2n}}B}\\
 \vdots & \vdots & \ddots & \vdots \\
{{a_{m1}}B}&{{a_{m2}}B}& \cdots &{{a_{mn}}B}
\end{array}} \right]\in {\mathbb{R}^{mp \times nq}}.\]
This definition implies that when $m=1$ and $q=1$, i.e, when $A$ is a row vector and $B$ is a column vector,
\begin{eqnarray}\label{1.8}
A \otimes B=BA.
\end{eqnarray}
From \cite[Chapter 4]{Horn91}, we have
\begin{eqnarray}
&&(A \otimes B)^{T}= (A^{T} \otimes B^{T}),\label{1.9}\\
&&{\rm vec}(AXB) = \left({B^{T}} \otimes A\right){\rm vec}(X),\label{1.10}\\
&&\Pi_{mn} {\rm vec}(A) = {\rm vec}({A^{T}}),\label{1.11}\\
&&\Pi_{pm} (A \otimes B) \Pi_{nq}= (B \otimes A),\nonumber
\end{eqnarray}
where $X\in {\mathbb{R}^{n \times p}}$, and $\Pi_{mn}$ is the {\em vec-permutation matrix} depending only on the orders $m$ and $n$. Especially, when $n=1$, i.e., $A$ is a column vector, then $\Pi_{nq}=I_q$ and hence
\begin{eqnarray}\label{1.12}
&&\Pi_{pm} (A \otimes B)= (B \otimes A).
\end{eqnarray}
In addition, the following result is also from \cite[Chapter 4]{Horn91}
\begin{eqnarray}\label{1.13}
(A \otimes B)(C \otimes D)=(AC) \otimes (BD),
\end{eqnarray}
where the matrices $C$ and $D$ are of suitable orders.

For the matrix pair $A ,B $ in \eqref{1.1} and \eqref{1.2}, there exist orthogonal matrices $U\in {\mathbb{R}^{m \times m}}$ and $V\in {\mathbb{R}^{s \times s}}$, and a nonsingular matrix $X\in {\mathbb{R}^{n \times n}}$ such that
\begin{eqnarray}\label{1.14}
A = U\Sigma X^{ - 1} ,\quad B = V\Lambda X^{ - 1},
\end{eqnarray}
where
\[
\Sigma  = \left[ {\begin{array}{*{20}c}
   {\Sigma _1 } & 0  \\
   0 & 0  \\
\end{array}} \right] = \left[ {\begin{array}{*{20}c}
   {I_{n - s} } & 0 & 0  \\
   0 & {S_A } & 0  \\
   0 & 0 & 0  \\
\end{array}} \right],\quad \Lambda  = \left[ {\begin{array}{*{20}c}
   0 & {\Lambda _1 }  \\
\end{array}} \right] = \left[ {\begin{array}{*{20}c}
   0 & {S_B } & {}  \\
   0 & {} & {I_{s - t} }  \\
\end{array}} \right]
\]
with
\[
S_A  = {\rm diag}(\alpha _1 , \cdots ,\alpha _t ),\ S_B  = {\rm diag}(\beta _1 , \cdots ,\beta _t ),\ \alpha _i ,\beta _i  > 0,\ \alpha _i^2  + \beta _i^2  = 1,
\]
and $t={\rm rank}(A)+s-n$. This decomposition is called the generalized SVD of a matrix pair (see e.g., \cite[p. 309]{Golub1}, \cite{Van1}). When $B=0$,  the generalized SVD can reduce to the SVD of $A$:
\begin{eqnarray}\label{1.15}
A = U\Sigma X^{ T} ,
\end{eqnarray}
where $X$ is orthogonal and $\Sigma _1 ={\rm diag}(\sigma _1 , \cdots ,\sigma _{{\rm rank}(A)})$ with $\sigma_i$ being the $i$-th singular value of $A$.

\section{The partial condition number}
Let $L\in {\mathbb{R}^{n \times k}}$ with $k\leq n$. We consider the following function
\begin{eqnarray*}
 g:\mathbb{R}^{m \times n}  \times \mathbb{R}^{s \times n}  \times \mathbb{R}^m  \times \mathbb{R}^s &\rightarrow&R^k\\
   (A,B,b,d)&\rightarrow&g(A,B,b,d) = L^{T} x(A,B,b,d) =  L^{T} (AP)^{\dag}  b+L^{T} B_A^{\dag}  d.
\end{eqnarray*}
From \cite{Cox, Li14}, it can be seen that the function $g$ is continuously Fr\'{e}chet differentiable in a neighborhood of $(A,B,b,d)$. Thus, denoting by $g{'}$ the Fr\'{e}chet derivative of $g$, and using the chain rules of composition of derivatives or from \cite{Cox,Li14}, we have
\begin{align*}
 &g{'}(A,B,b,d):\mathbb{R}^{m \times n}  \times \mathbb{R}^{s \times n}  \times \mathbb{R}^m  \times \mathbb{R}^s  \rightarrow R^k  \\
 &\quad\quad\quad\quad\quad\quad\quad\quad\  \ \  (\Delta A,\Delta B,\Delta b,\Delta d)\rightarrow g{'} (A,B,b,d){\circ}{(\Delta A,\Delta B,\Delta b,\Delta d)} \\
 &\quad\quad\quad\quad\quad\quad\quad\quad \quad\quad\quad\quad\quad\quad\quad\  \  = L^{T} ((AP)^{T} AP)^{\dag}  (\Delta A)^{T} r - L^{T} (AP)^{\dag}  (\Delta A)x + L^{T} (AP)^{\dag}  (\Delta b)\\
 &\quad\quad\quad\quad\quad\quad\quad\quad \quad\quad\quad\quad\quad\quad\quad\  \   \quad- L^{T} ((AP)^{T} AP)^{\dag}  (\Delta B)^{T} (AB_A^{\dag}  )^{T} r - L^{T} B_A^{\dag}  (\Delta B)x + L^{T} B_A^{\dag}  (\Delta d).
\end{align*}
Here, $ g{'} (A,B,b,d){\circ}{(\Delta A,\Delta B,\Delta b,\Delta d)}$ denotes that we apply the linear function $g{'}(A,B,b,d)$ to the perturbation variable  $(\Delta A,\Delta B,\Delta b,\Delta d)$ at the point $(A,B,b,d)$ and $r=b-Ax$ is called the residual vector.
Thus, according to \cite{Geu, Rice}, we have the absolute normwise condition number of $g$ at the point $(A,B,b,d)$ based on the weighted Frobenius norm \eqref{1.5}:
\begin{eqnarray}\label{2.1111}
\kappa_{LSE} (A,B,b,d) = \mathop {\max }\limits_{(\alpha_A\Delta A,\alpha_B\Delta B,\alpha_b\Delta b,\alpha_d\Delta d) \ne 0} \frac{{\left\| {g' (A,B,b,d){\circ} (\Delta A,\Delta B,\Delta b,\Delta d)} \right\|_2 }}{{\left\| {(\alpha_A\Delta A,\alpha_B\Delta B,\alpha_b\Delta b,\alpha_d\Delta d)} \right\|_F }}.
\end{eqnarray}
As mentioned in Section 1, the condition number $\kappa_{LSE} (A,B,b,d)$ is called the partial condition number of the LSE problem \eqref{1.1} with respect to $L$.

In the following, we provide an expression of $\kappa_{LSE} (A,B,b,d)$.
\begin{theorem}
The partial condition number of the LSE problem \eqref{1.1} with respect to $L$ is
\begin{eqnarray}\label{2.1}
\kappa_{LSE} (A,B,b,d) = \left\| {M_{g' } } \right\|_2 ,
\end{eqnarray}
where
\begin{eqnarray}\label{2.2}
M_{g' }  =\left[ {M_1 ,M_2 ,M_3 ,M_4 } \right]
\end{eqnarray}
with
\begin{eqnarray*}
 M_1  &=& \frac{{\left( {r^{T}  \otimes (L^{T} ((AP)^{T} AP)^{\dag}  )} \right)\Pi _{mn}  - {x^{T}  \otimes (L^{T} (AP)^{\dag}  )} }}{{\alpha _A }}, \\
 M_2  &=&  - \frac{{\left( {(r^{T} AB_A^{\dag}  ) \otimes (L^{T} ((AP)^{T} AP)^{\dag}  )} \right)\Pi _{sm}  +  {x^{T}  \otimes (L^{T} B_A^{\dag}  )} }}{{\alpha _B }}, \\
 M_3  &=& \frac{{L^{T} (AP)^{\dag}  }}{{\alpha _b }},\quad M_4  = \frac{{L^{T} B_A^{\dag}  }}{{\alpha _d }}.
\end{eqnarray*}
\end{theorem}
\emph{Proof}.
Applying the operator vec to $g{'} (A,B,b,d){\circ}{(\Delta A,\Delta B,\Delta b,\Delta d)}$ and using \eqref{1.10} and \eqref{1.11}, we have
\begin{align*}
 &g' (A,B,b,d){\circ}{(\Delta A,\Delta B,\Delta b,\Delta d)} ={\rm vec}(g' (A,B,b,d){\circ}{(\Delta A,\Delta B,\Delta b,\Delta d)}) \nonumber\\
 &\quad\quad\quad\quad\quad\quad\quad\quad\quad\quad\quad\quad\ = \left( {r^{T}  \otimes (L^{T} ((AP)^{T} AP)^{\dag}  )} \right)\Pi _{mn}{\rm vec}(\Delta A) - \left( {x^{T}  \otimes (L^{T} (AP)^{\dag}  )} \right){\rm vec}(\Delta A) \\
 &\quad\quad\quad\quad\quad\quad\quad\quad\quad\quad\quad\quad\   - \left( {(r^{T} AB_A^{\dag}  ) \otimes (L^{T} ((AP)^{T} AP)^{\dag}  )} \right)\Pi _{sm}{\rm vec}(\Delta B)-
  \left( {x^{T}  \otimes (L^{T} B_A^{\dag}  )} \right){\rm vec}(\Delta B)\\
 &\quad\quad\quad\quad\quad\quad\quad\quad\quad\quad\quad\quad\ +  L^{T} (AP)^{\dag}  (\Delta b) + L^{T} B_A^{\dag}  (\Delta d)\nonumber \\
 &\quad\quad\quad\quad\quad\quad\quad\quad\quad\quad\quad\quad\  = M_{g' } \left[ {\begin{array}{*{20}c}
   {\alpha_A {\rm vec}(\Delta A)}  \\
   {\alpha_B {\rm vec}(\Delta B)}  \\
   {\alpha_b (\Delta b)}  \\
   {\alpha_d (\Delta d)}  \\
\end{array}} \right].
\end{align*}
Thus, considering \eqref{1.5} and the fact that for any matrix $C$, $\left\|C\right\|_F=\left\|{\rm vec}(C)\right\|_2$,
\begin{eqnarray}\label{2.11111}
\kappa_{LSE} (A,B,b,d) = \mathop {\max }\limits_{(\alpha_A\Delta A,\alpha_B\Delta B,\alpha_b\Delta b,\alpha_d\Delta d) \ne 0} \frac{{\left\| {M_{g' } \left[ {\begin{array}{*{20}c}
   {\alpha_A {\rm vec}(\Delta A)}  \\
   {\alpha_B {\rm vec}(\Delta B)}  \\
   {\alpha_b (\Delta b)}  \\
   {\alpha_d (\Delta d)}  \\
\end{array}} \right]} \right\|_2 }}{{\left\| {\left[ {\begin{array}{*{20}c}
   {\alpha_A {\rm vec}(\Delta A)}  \\
   {\alpha_B {\rm vec}(\Delta B)}  \\
   {\alpha_b (\Delta b)}  \\
   {\alpha_d (\Delta d)}  \\
\end{array}} \right]} \right\|_2 }} = \left\| {M_{g' } } \right\|_2.
\end{eqnarray}
$\square$

\begin{remark}{
Setting $L=I_n$ and $\alpha_A=\alpha_B=\alpha_b=\alpha_d=1$ in \eqref{2.1}, and using the property on the spectral norm that for the matrices $C$ and $D$ of suitable orders, $\left\|[C,D]\right\|_2\leq\left\|C\right\|_2+\left\|D\right\|_2$, we have
\begin{eqnarray*}
\kappa_{LSE} (A,B,b,d) &\leq&   \left\|{\left( {r^{T}  \otimes  ((AP)^{T} AP)^{\dag}  } \right)\Pi _{mn}  -  {x^{T}  \otimes (AP)^{\dag}  } }\right\|_2\\
&\quad +& \left\|{\left( {(r^{T} AB_A^ \dag ) \otimes ((AP)^{T} AP)^{\dag}  } \right)\Pi _{sm}  +  {x^{T}  \otimes B_A^{\dag}  } }\right\|_2 +\left\|{(AP)^{\dag}  }\right\|_2+ \left\|{ B_A^{\dag}  }\right\|_2,
\end{eqnarray*}
which is essentially the same as the upper bound for the normwise condition number of the LSE problem \eqref{1.1} obtained in \cite{Cox}.
%, the differences stemming from different assumptions on the perturbations.
}
\end{remark}

Note that the expression of the partial condition number $\kappa_{LSE} (A,B,b,d)$ given in Theorem 2.1 contains Kronecker product. This introduces some large sparse matrices. The following theorem provides a closed formula of $\kappa_{LSE} (A,B,b,d)$ without Kronecker product.
\begin{theorem}
A closed formula of the partial condition number $\kappa_{LSE} (A,B,b,d)$ is given by
\begin{align}\label{2.3}
\kappa_{LSE} (A,B,b,d) = \left\|C\right\|_2^{1/2},
\end{align}
where
\begin{eqnarray}\label{2.4}
  C&=& \left( {\frac{{\left\| r \right\|_2^2 }}{{\alpha _A^2 }} + \frac{{\left\| {r^{T} AB_A^{\dag}  } \right\|_2^2 }}{{\alpha _B^2 }}} \right)L^{T} ((AP)^{T} AP)^{\dag}  )^2 L + \left( {\frac{{\left\| x \right\|_2^2 }}{{\alpha _A^2 }} + \frac{1}{{\alpha _b^2 }}} \right)L^{T} ((AP)^{T} AP)^{\dag}  )L \nonumber\\
  &\quad+& \left( {\frac{{\left\| x \right\|_2^2 }}{{\alpha _B^2 }} + \frac{1}{{\alpha _d^2 }}} \right)L^{T} B_A^{\dag}  (B_A^{\dag}  )^{T} L   + \frac{1}{{\alpha _B^2 }}L^{T} ((AP)^{T} AP)^{\dag}  xr^{T} AB_A^{\dag}  (B_A^{\dag}  )^{T} L \nonumber\\
  &\quad+& \frac{1}{{\alpha _B^2 }}L^{T} B_A^{\dag}  (B_A^{\dag}  )^{T}A^{T} rx^{T} ((AP)^{T} AP)^{\dag}  L .
\end{eqnarray}

\end{theorem}
\emph{Proof}. Noting
\begin{align*}
\left\| {M_{g' } } \right\|_2  = \left\| {M_{g' } M_{g' }^{T} } \right\|_2^{1/2}=\left\| M_1 M_1^{T}  + M_2 M_2^{T}  + M_3 M_3^{T}  + M_4 M_4^{T}\right\|_2^{1/2}
,
\end{align*}
and
\begin{align}\label{2.5}
M_3 M_3^{T}=
\frac{{L^{T} ((AP)^{T} AP)^{\dag}  )L}}{{\alpha _b^2 }}
, \quad M_4 M_4^{T}=\frac{{L^{T} B_A^{\dag}  (B_A^{\dag}  )^{T} L}}{{\alpha _d^2 }},
\end{align}
it suffices to obtain the expressions of $ M_1 M_1^{T}$ and $ M_2 M_2^{T}$.

Let
\begin{align*}
M_{11}  = \left( {r^{T}  \otimes (L^{T} ((AP)^{T} AP)^{\dag}  )} \right)\Pi _{mn} ,\quad M_{12}  =  {x^{T}  \otimes (L^{T} (AP)^{\dag}  )} .
\end{align*}
Then
\begin{align*}
M_1 M_1^{T}  = \frac{1}{{\alpha _A^2 }}\left( {M_{11} M_{11}^{T}  + M_{12} M_{12}^{T}  - M_{11} M_{12}^{T}  - M_{12} M_{11}^{T} } \right).
\end{align*}
By \eqref{1.9} and \eqref{1.13}, we have
\[
\begin{array}{l}
 M_{11} M_{11}^{T}  = \left( {r^{T}  \otimes (L^{T} ((AP)^{T} AP)^{\dag}  )} \right)\left( {r \otimes (((AP)^{T} AP)^{\dag}  L)} \right) = \left\| r \right\|_2^2 L^{T} ((AP)^{T} AP)^{\dag}  )^2 L, \\
 M_{12} M_{12}^{T}  = \left( {x^{T}  \otimes (L^{T} (AP)^{\dag}  )} \right)\left( {x \otimes (((AP)^{\dag}  )^{T} L)} \right) = \left\| x \right\|_2^2 L^{T} ((AP)^{T} AP)^{\dag}  )L.
 \end{array}
\]
Note that
\begin{eqnarray*}
 (AP)^{\dag}  r &=& (AP)^{\dag}  (b - Ax) = x - B_A^{\dag}  d - (AP)^{\dag}  Ax \quad\textrm{by \eqref{1.3}}\\
  &=&   x - B_A^{\dag}  Bx - (AP)^{\dag}  Ax=0 .
\end{eqnarray*}
The last equality in the above equation follows from the generalized SVD of the matrix pair $A,B$ in \eqref{1.14} and the expressions on $(AP)^{\dag}$ and $B_A^{\dag}$ in Remark \ref{rmk2.2} below. In fact,
\begin{eqnarray*}
B_A^{\dag}  B=X\Lambda^{\dag}\Lambda X^{-1}=X\left[ {\begin{array}{*{20}c}
   0 & 0   \\
   0 & I_s   \\
\end{array}} \right] X^{-1},\quad (AP)^{\dag}  A=X\left[ {\begin{array}{*{20}c}
   {I_{n - s} } & 0 & 0  \\
   0 & 0 & 0  \\
   0 & 0 & 0  \\
\end{array}} \right]\Sigma X^{-1}=X\left[ {\begin{array}{*{20}c}
   I_{n-s} & 0   \\
   0 & 0   \\
\end{array}} \right] X^{-1} ,
\end{eqnarray*}
which mean that $B_A^{\dag}  B+(AP)^{\dag}  A=I_n$ and hence $ x - B_A^{\dag}  Bx - (AP)^{\dag}  Ax=0$.
Thus, by \eqref{1.9}, \eqref{1.12}, and \eqref{1.13},
\begin{eqnarray*}
M_{11} M_{12}^{T}  &=& \left( {r^{T}  \otimes (L^{T} ((AP)^{T} AP)^{\dag}  )} \right)\left( {(((AP)^{\dag}  )^{T} L) \otimes x} \right)\quad\textrm{by \eqref{1.9} and \eqref{1.12}}\\
& =& (r^{T} ((AP)^{\dag}  )^{T} L) \otimes (L^{T} ((AP)^{T} AP)^{\dag}  x) \quad\textrm{by \eqref{1.13}}\\
&=&0=(
M_{12} M_{11}^{T})^{T}.
\end{eqnarray*}
As a result,
\begin{eqnarray}\label{2.6}
M_1 M_1^{T}  = \frac{1}{{\alpha _A^2 }}\left( {\left\| r \right\|_2^2 L^{T} ((AP)^{T} AP)^{\dag}  )^2 L + \left\| x \right\|_2^2 L^{T} ((AP)^{T} AP)^{\dag}  )L} \right).
\end{eqnarray}

Now, let
\[
M_{21}  = \left( {(r^{T} AB_A^{\dag}  ) \otimes (L^{T} ((AP)^{T} AP)^{\dag}  )} \right)\Pi _{sm} ,\quad M_{22}  = {x^{T}  \otimes (L^{T} B_A^{\dag}  )} .
\]
Then
\begin{eqnarray}\label{2.7}
M_2 M_2^{T}  = \frac{1}{{\alpha _B^2 }}\left( {M_{21} M_{21}^{T}  + M_{22} M_{22}^{T}  + M_{21} M_{22}^{T}  + M_{22} M_{21}^{T} } \right).
\end{eqnarray}
By \eqref{1.9} and \eqref{1.13}, we get
\begin{eqnarray}
 M_{21} M_{21}^{T} & =& \left( {(r^{T} AB_A^{\dag}  ) \otimes (L^{T} ((AP)^{T} AP)^{\dag}  )} \right)\left( {(r^{T} AB_A^{\dag}  )^{T}  \otimes (((AP)^{T} AP)^{\dag}  L)} \right)\nonumber \\
 &=& \left\| {r^{T} AB_A^{\dag}  } \right\|_2^2 L^{T} ((AP)^{T} AP)^{\dag}  )^2 L,\label{2.8} \\
 M_{22} M_{22}^{T}  &=& \left( {x^{T}  \otimes (L^{T} B_A^{\dag}  )} \right)\left( {x \otimes ((B_A^{\dag}  )^{T} L)} \right) = \left\| x \right\|_2^2 L^{T} B_A^{\dag}  (B_A^{\dag}  )^{T} L ,\label{2.9}
\end{eqnarray}
and by \eqref{1.9}, \eqref{1.12}, \eqref{1.13}, and \eqref{1.8}, we get
\begin{eqnarray}
 M_{21} M_{22}^{T}  &=& \left( {(r^{T} AB_A^{\dag}  ) \otimes (L^{T} ((AP)^{T} AP)^{\dag}  )} \right)\left( {((B_A^{\dag}  )^{T} L) \otimes x} \right)\quad\textrm{by \eqref{1.9} and \eqref{1.12}} \nonumber\\
  &=& (r^{T} AB_A^{\dag}  (B_A^{\dag}  )^{T} L) \otimes (L^{T} ((AP)^{T} AP)^{\dag}  x) \quad\textrm{by \eqref{1.13}}\nonumber\\
  &=& L^{T} ((AP)^{T} AP)^{\dag}  xr^{T} AB_A^{\dag}  (B_A^{\dag}  )^{T} L \label{2.10} \quad\textrm{by \eqref{1.8}}\\
 &=& (M_{22} M_{21}^{T})^{T}.\nonumber
  \end{eqnarray}
Substituting \eqref{2.8}--\eqref{2.10} into \eqref{2.7} gives
\begin{eqnarray}
M_2 M_2^{T}  &=&  \frac{1}{{\alpha _B^2 }}\left( {\left\| {r^{T} AB_A^{\dag}  } \right\|_2^2 L^{T} ((AP)^{T} AP)^{\dag}  )^2 L + \left\| x \right\|_2^2 L^{T} B_A^{\dag}  (B_A^{\dag}  )^{T} L} \right) \nonumber\\
  &\quad+& \frac{1}{{\alpha _B^2 }}L^{T} ((AP)^{T} AP)^{\dag}  xr^{T} AB_A^{\dag}  (B_A^{\dag}  )^{T} L + \frac{1}{{\alpha _B^2 }}L^{T} B_A^{\dag}  (B_A^{\dag}  )^{T} A^{T}rx^{T} ((AP)^{T} AP)^{\dag}  L.\label{2.11}
\end{eqnarray}
From \eqref{2.5}, \eqref{2.6}, and \eqref{2.11}, we have the desired result \eqref{2.3}. \quad $\square$

\begin{remark}{
When $B=0$ and $d=0$, that is, when the LSE problem reduces to the LLS problem \eqref{1.100}, $P=I_n$ and ${\rm rank}(A)=n$. Thus,
\[
((AP)^{T}(AP))^{\dag}   = (A^{T}A)^{-1},\quad B_A^{\dag}   = 0,
\]
and hence
\begin{eqnarray}\label{2.12}
  \kappa_{LLS} (A,b)&=&  \left\|\frac{{\left\| r \right\|_2^2 }}{{\alpha _A^2 }} L^{T} (A^{T}A)^{-2} L + \left( {\frac{{\left\| x \right\|_2^2 }}{{\alpha _A^2 }} + \frac{1}{{\alpha _b^2 }}} \right)L^{T} (A^{T}A)^{-1}L\right\|_2^{1/2} ,
\end{eqnarray}
which is the closed formula of the partial condition number of the LLS problem.

Furthermore, if $L$ is a column vector, i.e., $k=1$, then
\begin{eqnarray}
  \kappa_{LLS} (A,b)&=&  \left(\frac{{\left\| r \right\|_2^2 }}{{\alpha _A^2 }} L^{T} (A^{T}A)^{-2} L + \left( {\frac{{\left\| x \right\|_2^2 }}{{\alpha _A^2 }} + \frac{1}{{\alpha _b^2 }}} \right)L^{T} (A^{T}A)^{-1}L\right)^{1/2} \nonumber\\
  &=& \left(\frac{{\left\| r \right\|_2^2 }}{{\alpha _A^2 }} \left\|L^{T} (A^{T}A)^{-1} \right\|_2^2 + \left( {\frac{{\left\| x \right\|_2^2 }}{{\alpha _A^2 }} + \frac{1}{{\alpha _b^2 }}} \right)\left\|L^{T} A^{\dag}\right\|_2^2\right)^{1/2},\label{2.13}
\end{eqnarray}
which is just the result given in Corollary 1 in \cite{Ari}.
}
\end{remark}

\begin{remark}\label{rmk2.2}{
Using the generalized SVD of the matrix pair $A,B$ in \eqref{1.14}, and (3.3), (3.4), and (3.15) in \cite{Wei921}, we have
\begin{eqnarray*}
(AP)^{\dag}   &=& X(\Sigma (I_n  - \Lambda ^{\dag}  \Lambda ))^{\dag}  U^{T}\nonumber\\
  &=& X\left( {\left[ {\begin{array}{*{20}c}
   {I_{n - s} } & 0 & 0  \\
   0 & {S_A } & 0  \\
   0 & 0 & 0  \\
\end{array}} \right]\left( {I_n  - \left[ {\begin{array}{*{20}c}
   0 & 0  \\
   {S_B^{ - 1} } & 0  \\
   0 & {I_{s - t} }  \\
\end{array}} \right]\left[ {\begin{array}{*{20}c}
   0 & {S_B } & 0  \\
   0 & 0 & {I_{s - t} }  \\
\end{array}} \right]} \right)} \right)^{\dag}  U^{T} \nonumber \\
  &=& X\left( {\left[ {\begin{array}{*{20}c}
   {I_{n - s} } & 0 & 0  \\
   0 & {S_A } & 0  \\
   0 & 0 & 0  \\
\end{array}} \right]\left[ {\begin{array}{*{20}c}
   {I_{n - s} } & 0 & 0  \\
   0 & 0 & 0  \\
   0 & 0 & 0  \\
\end{array}} \right]} \right)^{\dag}  U^{T}  = X\left[ {\begin{array}{*{20}c}
   {I_{n - s} } & 0 & 0  \\
   0 & 0 & 0  \\
   0 & 0 & 0  \\
\end{array}} \right]U^{T}
\end{eqnarray*}
and
\begin{eqnarray*}
B_A^{\dag}  & =& (I_n  - (AP)^{\dag}  A)X\Lambda ^{\dag}  V^{T}\nonumber\\
  &= &\left( {I_n  - X\left[ {\begin{array}{*{20}c}
   {I_{n - s} } & 0 & 0  \\
   0 & 0 & 0  \\
   0 & 0 & 0  \\
\end{array}} \right]U^{T} U\left[ {\begin{array}{*{20}c}
   {I_{n - s} } & 0 & 0  \\
   0 & {S_A } & 0  \\
   0 & 0 & 0  \\
\end{array}} \right]X^{ - 1} } \right)X\Lambda ^{\dag}  V^{T} \nonumber \\
  &=&  X\left[ {\begin{array}{*{20}c}
   0 & 0 & 0  \\
   0 & {I_t } & 0  \\
   0 & 0 & {I_{s - t} }  \\
\end{array}} \right]\left[ {\begin{array}{*{20}c}
   0 & 0  \\
   {S_B^{ - 1} } & 0  \\
   0 & {I_{s - t} }  \\
\end{array}} \right]V^{T}  = X\left[ {\begin{array}{*{20}c}
   0 & 0  \\
   {S_B^{ - 1} } & 0  \\
   0 & {I_{s - t} }  \\
\end{array}} \right]V^{T}  = X\Lambda ^{\dag}  V^{T}.
\end{eqnarray*}
Then
\begin{eqnarray}
&(AP)^{\dag}  \left( {(AP)^{\dag}  } \right)^{T} =X_1X_1^{T} , \quad AB_A^{\dag}=U\left[ {\begin{array}{*{20}c}
   0 & 0  \\
   {S_AS_B^{ - 1} } & 0  \\
   0 & 0  \\
\end{array}} \right]V^{T}=U_2\left[ {\begin{array}{*{20}c}
   {S_AS_B^{ - 1} } & 0  \\
   0 & 0  \\
\end{array}} \right]V^{T},\label{2.14}\\
& B_A^{\dag}  (B_A^{\dag}  )^{T}  = X\left[ {\begin{array}{*{20}c}
   0 & 0 & 0  \\
   0 & {S_B^{ - 2} } & 0  \\
   0 & 0 & {I_{s - t} }  \\
\end{array}} \right]X^{T}=X_2\left[ {\begin{array}{*{20}c}
     {S_B^{ - 2} } & 0  \\
    0 & {I_{s - t} }  \\
\end{array}} \right]X^{T}_2 , \label{2.15}\\
& AB_A^{\dag}  (B_A^{\dag}  )^{T}   = U\left[ {\begin{array}{*{20}c}
   0 & 0 & 0  \\
   0 & {S_A S_B^{ - 2} } & 0  \\
   0 & 0 & 0  \\
\end{array}} \right]X^{T}= U_2\left[ {\begin{array}{*{20}c}
   {S_A S_B^{ - 2} } & 0  \\
   0 & 0  \\
\end{array}} \right]X^{T}_2,\label{2.16}
\end{eqnarray}
where $X=[X_1,X_2]$ with $ X_1\in \mathbb{R}^{n \times (n-s)}$ and $ X_2\in \mathbb{R}^{n \times s}$, and $U=[U_1,U_2]$ with $ U_1\in \mathbb{R}^{m \times (n-s)}$ and $ U_2\in \mathbb{R}^{m \times (m-n+s)}$. Substituting \eqref{2.14}--\eqref{2.16} into \eqref{2.4} yields
\begin{align}
 &C = \left( {\frac{{\left\| r \right\|_2^2 }}{{\alpha _A^2 }} + \frac{{\left\| {r^{T} U_2\left[ {\begin{array}{*{20}c}
   {S_A S_B^{ - 1} } & 0  \\
   0 & 0  \\
\end{array}} \right] } \right\|_2^2 }}{{\alpha _B^2 }}} \right)L^{T} (X_1X_1^{T})^2 L + L^{T} (X_1S_1X_1^{T}+X_2S_2X_2^{T}) L  \nonumber \\
&\quad   + \frac{1}{{\alpha _B^2 }}L^{T} X_1X^{T}_1 xr^{T} U_2\left[ {\begin{array}{*{20}c}
   {S_A S_B^{ - 2} } & 0  \\
   0 & 0  \\
\end{array}} \right]X^{T}_2 L +\frac{1}{{\alpha _B^2 }}L^{T} X_2\left[ {\begin{array}{*{20}c}
   {S_A S_B^{ - 2} } & 0  \\
   0 & 0  \\
\end{array}} \right]U^{T}_2 rx^{T} X_1X^{T}_1 L\label{2.17}
\end{align}
with
\begin{eqnarray}\label{2.18}
S_1=\left( {\frac{{\left\| x \right\|_2^2 }}{{\alpha _A^2 }} + \frac{1}{{\alpha _b^2 }}} \right)I_{n-s},\quad S_2=\left( {\frac{{\left\| x \right\|_2^2 }}{{\alpha _B^2 }} + \frac{1}{{\alpha _d^2 }}} \right){\Lambda _1^{ - 2} }.
\end{eqnarray}

In particular, when $B=0$, the generalized SVD \eqref{1.14} reduces to the SVD of $A$ \eqref{1.15}. In this case, $P=I_n$ and ${\rm rank}(A)=n$. Hence, we have
\[
(AP)^{\dag}   = X\Sigma ^{\dag}  U^{T} =X[\Sigma ^ {-1}_1,0]  U^{T},\quad B_A^{\dag}   = 0,
\]
and 
\[
(AP)^{\dag}  \left( {(AP)^{\dag}  } \right)^{T}    = X\Sigma ^ {-2}_1 X^{T}.
\]
Thus,
 \[
C = \frac{{\left\| r \right\|_2^2 }}{{\alpha _A^2 }}L^{T}X \Sigma ^ {-4}_1 X^{T}L + \left( {\frac{{\left\| x \right\|_2^2 }}{{\alpha _A^2 }} + \frac{1}{{\alpha _b^2 }}} \right)L^{T} X\Sigma ^ {-2}_1X^{T} L.
\]
As a result, we get a closed formula of the partial condition number of the LLS problem based on the SVD of $A$:
\begin{eqnarray}\label{2.19}
\kappa_{LLS} (A,b)  = \left\| {SX^{T} L} \right\|_2,
\end{eqnarray}
where $S$ is a diagonal matrix with the $i$-th diagonal element being
\[
S_{ii}  =\frac{1}{\sigma _i} \sqrt {\frac{{\sigma _i^{-2}\left\| r \right\|_2^2+\left\| x \right\|_2^2 }}{{\alpha_A ^2 }} + \frac{1}{{\alpha_b ^2 }}} .
\]
The closed formula \eqref{2.19} is just the one given in Theorem 1 in \cite{Ari}, where it was derived by a different approach.}
\end{remark}

\section{The structured partial condition number}
Suppose that $\mathbb{S}_1\subseteq\mathbb{R}^{m\times n}$ and $\mathbb{S}_2\subseteq \mathbb{R}^{s\times n}$ are two linear subspaces, which consist of two classes of structured matrices, respectively. From \cite{Hig,Lijia,Rump1}, we have that if $A\in \mathbb{S}_1$ and $B \in \mathbb{S}_2$, then
\begin{eqnarray}\label{3.1}
{\rm vec}(A)=\Phi_{\mathbb{S}_1}s_1,\quad {\rm vec}(B)=\Phi_{\mathbb{S}_2}s_2,
\end{eqnarray}
where $\Phi_{\mathbb{S}_1}\in \mathbb{R}^{mn\times k_1}$ and $\Phi_{\mathbb{S}_2}\in \mathbb{R}^{sn\times k_2}$ are the fixed structure matrices reflecting the structures of $\mathbb{S}_1$ and $\mathbb{S}_2$, respectively, and $s_1\in \mathbb{R}^{k_1}$ and $s_2\in \mathbb{R}^{k_2}$ are the vectors of the independent parameters in the structured matrices, respectively. Based on the above explanation, the structured perturbations $\Delta A\in \mathbb{S}_1$ and $\Delta B \in \mathbb{S}_2$ can be written as
\begin{eqnarray}\label{3.2}
{\rm vec}(\Delta A)=\Phi_{\mathbb{S}_1}(\Delta s_1),\quad {\rm vec}(\Delta B)=\Phi_{\mathbb{S}_2}(\Delta s_2),
\end{eqnarray}
where $\Delta s_1\in \mathbb{R}^{k_1}$ and $\Delta s_2\in \mathbb{R}^{k_2}$ can be regarded as the perturbations of $s_1$ and $s_2$, respectively.

Now we present the definition of the structured partial condition number of the LSE problem \eqref{1.1}:
\begin{eqnarray*}
\kappa_{LSE}^S (A,B,b,d) = \mathop {\max }\limits_{(\alpha_A \Delta A,\alpha_B \Delta B,\alpha_b \Delta b,\alpha_d \Delta d) \ne 0\hfill \atop
  \scriptstyle
\Delta A\in \mathbb{S}_1, \Delta B\in \mathbb{S}_2\hfill} \frac{{\left\| {g' (A,B,b,d){\circ} (\Delta A,\Delta B,\Delta b,\Delta d)} \right\|_2 }}{{\left\| {(\alpha_A\Delta A,\alpha_B\Delta B,\alpha_b\Delta b,\alpha_d\Delta d)} \right\|_F }},
\end{eqnarray*}
which is a natural variant of the partial condition number in \eqref{2.1111}. From \eqref{2.11111}, it follows that
\begin{eqnarray}\label{3.3}
\kappa_{LSE}^S (A,B,b,d) = \mathop {\max }\limits_{(\alpha_A \Delta A,\alpha_B \Delta B,\alpha_b \Delta b,\alpha_d \Delta d) \ne 0\hfill \atop
  \scriptstyle
\Delta A\in \mathbb{S}_1, \Delta B\in \mathbb{S}_2\hfill}\frac{{\left\| {M_{g' } \left[ {\begin{array}{*{20}c}
   {\alpha_A {\rm vec}(\Delta A)}  \\
   {\alpha_B {\rm vec}(\Delta B)}  \\
   {\alpha_b (\Delta b)}  \\
   {\alpha_d (\Delta d)}  \\
\end{array}} \right]} \right\|_2 }}{{\left\| {\left[ {\begin{array}{*{20}c}
   {\alpha_A {\rm vec}(\Delta A)}  \\
   {\alpha_B {\rm vec}(\Delta B)}  \\
   {\alpha_b (\Delta b)}  \\
   {\alpha_d (\Delta d)}  \\
\end{array}} \right]} \right\|_2 }} .
\end{eqnarray}
Considering \eqref{3.2}, we have
\begin{eqnarray*}
 \left[ {\begin{array}{*{20}c}
   { {\rm vec}(\Delta A)}  \\
   { {\rm vec}(\Delta B)}  \\
   { \Delta b}  \\
     { \Delta d}  \\
\end{array}} \right]=\left[ {\begin{array}{*{20}c}
   \Phi_{\mathbb{S}_1} &0 &0&0 \\
   0&\Phi_{\mathbb{S}_2}&0&0  \\
   0&0& I_m&0\\
   0&0&0&I_s\\
\end{array}} \right]\left[ {\begin{array}{*{20}c}
   \Delta s _1 \\
   \Delta s_2 \\
    { \Delta b}  \\
     { \Delta d}  \\
\end{array}} \right].
\end{eqnarray*}
Substituting the above equation into \eqref{3.3} yields
\begin{eqnarray}\label{3.4}
\kappa_{LSE}^S (A,B,b,d) = \mathop {\max }\limits_{(\alpha_A (\Delta s_1),{\alpha_B (\Delta s_2)},{\alpha_b (\Delta b)},{\alpha_d (\Delta d)} ) \ne 0}\frac{{\left\| {M_{g' } \left[ {\begin{array}{*{20}c}
   \Phi_{\mathbb{S}_1} &0 &0&0 \\
   0&\Phi_{\mathbb{S}_2}&0&0  \\
   0&0& I_m&0\\
   0&0&0&I_s\\
\end{array}} \right]\left[ {\begin{array}{*{20}c}
   {\alpha_A (\Delta s_1)}  \\
   {\alpha_B (\Delta s_2)}  \\
   {\alpha_b (\Delta b)}  \\
   {\alpha_d (\Delta d)}  \\
\end{array}} \right]} \right\|_2 }}{{\left\| {\left[ {\begin{array}{*{20}c}
   \Phi_{\mathbb{S}_1} &0 &0&0 \\
   0&\Phi_{\mathbb{S}_2}&0&0  \\
   0&0& I_m&0\\
   0&0&0&I_s\\
\end{array}} \right]\left[ {\begin{array}{*{20}c}
   {\alpha_A (\Delta s_1)}  \\
   {\alpha_B (\Delta s_2)}  \\
   {\alpha_b (\Delta b)}  \\
   {\alpha_d (\Delta d)}  \\
\end{array}} \right]} \right\|_2 }} .
\end{eqnarray}
Note that
\begin{eqnarray*}
{{\left\| {\left[ {\begin{array}{*{20}c}
   \Phi_{\mathbb{S}_1} &0 &0&0 \\
   0&\Phi_{\mathbb{S}_2}&0&0  \\
   0&0& I_m&0\\
   0&0&0&I_s\\
\end{array}} \right]\left[ {\begin{array}{*{20}c}
   {\alpha_A (\Delta s_1)}  \\
   {\alpha_B (\Delta s_2)}  \\
   {\alpha_b (\Delta b)}  \\
   {\alpha_d (\Delta d)}  \\
\end{array}} \right]} \right\|_2 }}=\left\|\left[ {\begin{array}{*{20}c}
   {\alpha_A (\Delta s_1)}  \\
   {\alpha_B (\Delta s_2)}  \\
   {\alpha_b (\Delta b)}  \\
   {\alpha_d (\Delta d)}  \\
\end{array}} \right]^{T} {\left[ {\begin{array}{*{20}c}
   \Phi_{\mathbb{S}_1}^{T}\Phi_{\mathbb{S}_1} &0 &0&0 \\
   0&\Phi_{\mathbb{S}_2}^{T}\Phi_{\mathbb{S}_2} &0&0 \\
   0&0& I_m&0\\
   0&0&0&I_s\\
\end{array}} \right]\left[ {\begin{array}{*{20}c}
   {\alpha_A (\Delta s_1)}  \\
   {\alpha_B (\Delta s_2)}  \\
   {\alpha_b (\Delta b)}  \\
   {\alpha_d (\Delta d)}  \\
\end{array}} \right]} \right\|_2^{1/2}
\end{eqnarray*}
and the structured matrices $\Phi_{\mathbb{S}_1}$ and $\Phi_{\mathbb{S}_2}$ are column orthogonal \cite{Lijia}. Then
\begin{eqnarray}\label{3.5}
{{\left\| {\left[ {\begin{array}{*{20}c}
   \Phi_{\mathbb{S}_1} &0 &0&0 \\
   0&\Phi_{\mathbb{S}_2}&0&0  \\
   0&0& I_m&0\\
   0&0&0&I_s\\
\end{array}} \right]\left[ {\begin{array}{*{20}c}
   {\alpha_A (\Delta s_1)}  \\
   {\alpha_B (\Delta s_2)}  \\
   {\alpha_b (\Delta b)}  \\
   {\alpha_d (\Delta d)}  \\
\end{array}} \right]} \right\|_2 }}=\left\| {\left[ {\begin{array}{*{20}c}
   D_1 &0 &0&0 \\
   0&D_2&0&0  \\
   0&0& I_m&0\\
   0&0&0&I_s\\
\end{array}} \right]\left[ {\begin{array}{*{20}c}
   {\alpha_A (\Delta s_1)}  \\
   {\alpha_B (\Delta s_2)}  \\
   {\alpha_b (\Delta b)}  \\
   {\alpha_d (\Delta d)}  \\
\end{array}} \right]} \right\|_2,
\end{eqnarray}
where $D_1={\rm diag}(w_1)$ and $D_2={\rm diag}(w_2)$ with
\begin{eqnarray*}
w_1=\left[\left\|\Phi_{\mathbb{S}_1}(1,:)\right\|_2,\cdots,\left\|\Phi_{\mathbb{S}_1}(k_1,:)\right\|_2\right],\quad w_2=\left[\left\|\Phi_{\mathbb{S}_2}(1,:)\right\|_2,\cdots,\left\|\Phi_{\mathbb{S}_2}(k_2,:)\right\|_2\right].
\end{eqnarray*}
Here, the Matlab notation is used.
Combining \eqref{3.4} and \eqref{3.5} implies
\begin{align*}
&\kappa_{LSE}^S (A,B,b,d)\\
& = \mathop {\max }\limits_{(\alpha_A (\Delta s_1),{\alpha_B (\Delta s_2)},{\alpha_b (\Delta b)},{\alpha_d (\Delta d)} ) \ne 0}\frac{{\left\| {M_{g' } \left[ {\begin{array}{*{20}c}
   \Phi_{\mathbb{S}_1}D_1^{-1} &0 &0&0 \\
   0&\Phi_{\mathbb{S}_2}D_2^{-1}&0&0  \\
   0&0& I_m&0\\
   0&0&0&I_s\\
\end{array}} \right]\left[ {\begin{array}{*{20}c}
   D_1 &0 &0&0 \\
   0&D_2&0&0  \\
   0&0& I_m&0\\
   0&0&0&I_s\\
\end{array}} \right]\left[ {\begin{array}{*{20}c}
   {\alpha_A (\Delta s_1)}  \\
   {\alpha_B (\Delta s_2)}  \\
   {\alpha_b (\Delta b)}  \\
   {\alpha_d (\Delta d)}  \\
\end{array}} \right]} \right\|_2 }}{{\left\| {\left[ {\begin{array}{*{20}c}
   D_1 &0 &0&0 \\
   0&D_2&0&0  \\
   0&0& I_m&0\\
   0&0&0&I_s\\
\end{array}} \right]\left[ {\begin{array}{*{20}c}
   {\alpha_A (\Delta s_1)}  \\
   {\alpha_B (\Delta s_2)}  \\
   {\alpha_b (\Delta b)}  \\
   {\alpha_d (\Delta d)}  \\
\end{array}} \right]} \right\|_2 }} .
\end{align*}
Then we can derive the expression of the structured partial condition number of the LSE problem, which is presented in the following theorem.
\begin{theorem}
The structured partial condition number of the LSE problem \eqref{1.1} with respect to $L$ and the structures $\mathbb{S}_1$ and $\mathbb{S}_2$ is
\begin{eqnarray}\label{3.6}
\kappa_{LSE}^S (A,B,b,d) = \left\| {M_{g' } } \left[ {\begin{array}{*{20}c}
   \Phi_{\mathbb{S}_1}D_1^{-1} &0 &0&0 \\
   0&\Phi_{\mathbb{S}_2}D_2^{-1}&0&0  \\
   0&0& I_m&0\\
   0&0&0&I_s\\
\end{array}} \right]\right\|_2 ,
\end{eqnarray}
where $M_{g' }$ is given in \eqref{2.2}.
\end{theorem}
\begin{remark}\label{rmk3.2}
{
It is easy to verify that
$$ \left[ {\begin{array}{*{20}c}
   \Phi_{\mathbb{S}_1}D_1^{-1} &0 &0&0 \\
   0&\Phi_{\mathbb{S}_2}D_2^{-1}&0&0  \\
   0&0& I_m&0\\
   0&0&0&I_s\\
\end{array}} \right]$$
is column orthonormal. Thus,
$$\left\| {M_{g' } } \left[ {\begin{array}{*{20}c}
   \Phi_{\mathbb{S}_1}D_1^{-1} &0 &0&0 \\
   0&\Phi_{\mathbb{S}_2}D_2^{-1}&0&0  \\
   0&0& I_m&0\\
   0&0&0&I_s\\
\end{array}} \right]\right\|_2\leq \left\| {M_{g' } }\right\|_2.$$
That is, the structured partial condition number is always tighter than the unstructured one. This fact can also be seen from the definitions of these two condition numbers. As done in \cite{Rump1, Xu}, it is valuable to discuss the ratio between the structured and unstructured partial condition numbers of the
LSE problem in detail. We won't go that far in this paper, and only provide a numerical example in Section 5 to show that the structured partial condition number is indeed tighter than the unstructured one.
}
\end{remark}

\begin{remark}{
When $B=0$ and $d=0$, we have the structured partial condition number of the LLS problem and its upper bound:
\begin{eqnarray}
\kappa_{LLS}^S (A,b) &=& \left\|\left[ {\frac{{\left( {r^{T}  \otimes (L^{T} (A^{T}A)^{ - 1} )} \right)\Pi _{mn}  - x^{T}  \otimes (L^{T} A^{ \dag})}}{\alpha_A }\Phi_{\mathbb{S}_1}D_1^{-1} ,\frac{{L^{T} A^{ \dag} }}{\alpha_b }}  \right]\right\|_2 \label{3.7}\\
&\leq&\left\|\left[ {\frac{{\left( {r^{T}  \otimes (L^{T} (A^{T}A)^{ - 1} )} \right)\Pi _{mn}  - x^{T}  \otimes (L^{T} A^{ \dag})}}{\alpha_A } ,\frac{{L^{T} A^{ \dag} }}{\alpha_b }}  \right]\right\|_2 ,\label{3.8}
\end{eqnarray}
where the upper bound \eqref{3.8} is just the unstructured partial condition number of the LLS problem. Here, it should be pointed out that the structured condition number of the LLS problem derived from \eqref{3.7} by setting $L=I_n$ and $\alpha_A=\alpha_b=1$ is a little different from the ones in \cite{Xu} because two additional conditions are added besides the structure requirement in \cite{Xu}.
}
\end{remark}

\begin{remark}{ We only consider the linear structures of the matrices $A$ and $B$ in this section. Similarly, the linear structures of the vectors $b$ and $d$ can also be put into the partial condition number. Furthermore, inspired by \cite{Cucker,Lijia,Rump2}, exploring the structured mixed and componentwise condition numbers of the LSE problem will be interesting. We will investigate this problem in the future research.
}
\end{remark}

\section{Statistical condition estimates}
We first provide a statistical estimate of the partial condition number by using the probabilistic spectral norm estimator. This estimator was proposed by Hochstenbach \cite{Hochs13} and can estimate the spectral norm of a matrix reliably. In more detail, the analysis of the estimator in \cite{Hochs13} suggests that the spectral norm of a matrix can be contained in a small interval $[\alpha_1 , \alpha_2]$ with high probability, where $\alpha_1$ is the guaranteed lower bound of the spectral norm of the matrix derived by the famous Lanczos bibdiagonalization method \cite{Golub0} and $\alpha_2$ is the probabilistic upper bound with probability at least $1-\varepsilon$ with $\varepsilon\ll 1$ derived by finding the largest zero of a polynomial. Meanwhile, we can require $\alpha_2/ \alpha_1\leqslant 1+\delta$ with $\delta$ being a user-chosen parameter. Based on the above estimator, we can devise Algorithm \ref{AlgorithmPCE}.

\begin{algorithm}[htbp]{ 
\caption{Probabilistic spectral norm estimator for the partial condition number \eqref{2.3}}\label{AlgorithmPCE}
\begin{enumerate}
  \item Generate a starting vector $v_0$ from $\mathcal{U}(S_{q-1})$ with $q=n^2$. Hereafter, $\mathcal{U}(S_{q-1})$ denotes the uniform distribution over unit sphere $S_{q-1}$ in $R^{q}$.
  \item Compute the guaranteed lower bound $\alpha_1$ and the probabilistic upper bound $\alpha_2$ of $\left\|C\right\|_2$ by probabilistic spectral norm estimator, where $C$ is given in \eqref{2.4} or \eqref{2.17}.
  \item Estimate the partial condition number \eqref{2.3} by
  \begin{equation*}
  {\kappa}_{PLSE}(A,B,b,d)=\sqrt{\frac{\alpha_1+\alpha_2}{2}}.
  \end{equation*}
\end{enumerate}}
\end{algorithm}
\begin{remark}{
In the practical implementation of Algorithm \ref{AlgorithmPCE}, explicitly forming matrix $C$ is not necessary because what we really need is the product of a random vector with the matrix $C$ or $C^T$. Hence, some techniques in solving linear system can be employed to reduce the computational burden. Furthermore, it is worthy to point out that Algorithm \ref{AlgorithmPCE} is also applicable to estimating the partial structured condition number \eqref{3.6} since it is also the spectral norm of a matrix.
}
\end{remark}

Now we introduce an alternative approach based on the SSCE method \cite{Babo,Kenney94} for estimating the normwise condition number of the solution $x(A,B,b,d)$. Denote by $\kappa_{LSEi}(A,B,b,d)$ the normwise condition number of the function $z_i^{T}x(A,B,b,d)$, where $z_i$s are chosen from $\mathcal{U}(S_{n-1})$ and are orthogonal. Then, from \eqref{2.4}, we have
\begin{eqnarray}
  \kappa_{LSEi}^2 (A,B,b,d)&=& \left( {\frac{{\left\| r \right\|_2^2 }}{{\alpha _A^2 }} + \frac{{\left\| {r^{T} AB_A^{\dag}  } \right\|_2^2 }}{{\alpha _B^2 }}} \right)z_i^{T} ((AP)^{T} AP)^{\dag}  )^2 z_i + \left( {\frac{{\left\| x \right\|_2^2 }}{{\alpha _A^2 }} + \frac{1}{{\alpha _b^2 }}} \right)z_i^{T} ((AP)^{T} AP)^{\dag}  )z_i \nonumber\\
  &&+ \left( {\frac{{\left\| x \right\|_2^2 }}{{\alpha _B^2 }} + \frac{1}{{\alpha _d^2 }}} \right)z_i^{T} B_A^{\dag}  (B_A^{\dag}  )^{T} z_i   + \frac{2}{{\alpha _B^2 }}z_i^{T} ((AP)^{T} AP)^{\dag}  xr^{T} AB_A^{\dag}  (B_A^{\dag}  )^{T} z_i.\label{4.1}
\end{eqnarray}
The analysis based on SSCE method in \cite{Babo} shows that
\begin{eqnarray}\label{4.2}
{\kappa} _{SLSE}(A,B,b,d) = \frac{\omega_q}{\omega_n}\sqrt{\sum_{i=1}^q \kappa_{LSEi}^2 (A,B,b,d)}
\end{eqnarray}
is a good estimate of the normwise condition number of the LSE problem \eqref{1.1}. In the above expression, $\omega_q$ is the Wallis factor with $\omega_1=1$, $\omega_2={2}/{\pi}$, and
$$\omega_q  =\left\{
  \begin{array}{ll}
    \frac{1\cdot 3 \cdot5 \cdots (q-2)}{2\cdot 4\cdot 6\cdots (q-1)}, & \hbox{for $q$ odd,} \\
    \frac{2}{\pi}\frac{2\cdot 4\cdot 6\cdots (q-2)}{3\cdot 5\cdot 7\cdots (q-1)}, & \hbox{for $q$ even,}
  \end{array}
\right.\textrm{ when } q>2.$$
It can be approximated by
\begin{eqnarray}\label{4.3}
\omega_q\approx \sqrt{\frac{2}{\pi(q-\frac{1}{2})}}
\end{eqnarray}
with high accuracy. In summary, we can propose Algorithm \ref{AlgorithmSSCE}.

\begin{algorithm}[htbp]{
\caption{SSCE method for the normwise condition number of the LSE solution}\label{AlgorithmSSCE}
\begin{enumerate}
  \item Generate $q$ vectors $z_1,\cdots, z_q$ from $\mathcal{U}(S_{n-1})$, and orthonormalize these vectors using the QR facotization.
  \item For $i=1,\cdots, q$, compute $\kappa_{LSEi}^2 (A,B,b,d)$ by \eqref{4.1}.
  \item Approximate $\omega_q$ and $\omega_n$ by \eqref{4.3} and estimate the normwise condition number by \eqref{4.2}.
\end{enumerate}}
\end{algorithm}
\begin{remark}{ 
In Algorithm \ref{AlgorithmSSCE}, $\kappa_{LSEi}^2 (A,B,b,d)$ is computed by the equation \eqref{4.1}. In practice, the computation of  $\kappa_{LSEi}^2 (A,B,b,d)$ should rely on the intermediate results of the process for solving the LSE problem to reduce the computational burden. Just as carried out in \cite{Babo}, where the estimate is computed by using the $R$ factor of QR decomposition, it is better to compute $\kappa_{LSEi}^2 (A,B,b,d)$ through a formula descended from \eqref{2.17} instead of \eqref{2.4} if we solve the LSE problem by generalized SVD.
}\end{remark}

\section{Numerical experiments}
In this section, we will present two numerical examples to illustrate the reliability of the statistical condition estimates proposed in Section 4 and to compare the structured condition number and the  unstructured one, respectively. In these two examples, we will set $\alpha_A=\alpha_B=\alpha_b=\alpha_d=1$ and the matrix $L$ be the identity matrix.

\begin{example}\label{example1}
{
Similar to \cite{Paige82}, we generate the example as follows. Let $u_1\in \mathbb{R}^{m}$, $u_2\in \mathbb{R}^{s}$, and $v_1, v_2 \in \mathbb{R}^{n}$ be unit random vectors, and set
\begin{eqnarray*}
  A &=& U_1\begin{bmatrix}
             D_1 \\
             0 \\
           \end{bmatrix}V_1, \
           B=U_2\begin{bmatrix}
                  D_2 & 0 \\
                \end{bmatrix}V_2 ,\ U_i=I_{m(s)}-2u_iu_i^T,\ \mathrm{and} \ V_i=I_n-2v_iv_i^T,
\end{eqnarray*}
where $D_1=n^{-l_1}\mathrm{diag}(n^{l_1},(n-1)^{l_1},\cdots,1)$ and $D_2=s^{-l_2}\mathrm{diag}(s^{l_2},(s-1)^{l_2},\cdots,1)$. Let the solution $x$ be $x=(1, 2^2,\cdots, n^2)^T$ and the residual vector $r=b-Ax$ be a random vector of specified norm. Thus, letting $b=Ax+r$ and $d=Bx$ gives the desired LSE problem, and it is easy to check that the condition numbers of $A$ and $B$ are $\kappa(A)=n^{l_1}$ and $\kappa(B)=s^{l_2}$, respectively. Recall that for any matrix $C$, its condition number $\kappa(C)$ is defined by $\kappa(C)=\left\|C\right\|_2\left\|C^\dag\right\|_2$.

In our numerical experiments, we set $m=100$, $n=80$ and $s=50$, and choose the parameters $\epsilon=0.001$, $\delta= 0.01$ in Algorithm \ref{AlgorithmPCE} and $q=2$ in Algorithm \ref{AlgorithmSSCE}.  By varying the condition numbers of $A$ and $ B$, and the residual's norm $\|r\|_2$, we test the performance of Algorithms \ref{AlgorithmPCE} and \ref{AlgorithmSSCE}. More precisely, for each pair of $\kappa(A)$ and $\|r\|_2$  with a fixed $\kappa(B)$, $500$ random LSE problems are generated and used for the test. The numerical results on mean and variance of the ratios between the statistical condition estimate and the exact condition number defined as $$r_{ssce}={\kappa} _{SLSE}(A,B,b,d)/{\kappa} _{LSE}(A,B,b,d),\ r_{pce}={\kappa} _{PLSE}(A,B,b,d)/{\kappa} _{LSE}(A,B,b,d)$$ are reported in Tables \ref{Table1}.

From Table \ref{Table1}, one can easily find that in general both Algorithms \ref{AlgorithmPCE} and \ref{AlgorithmSSCE} can give reliable estimates of the normwise condition number. In comparison, Algorithm \ref{AlgorithmPCE} performs more stable since the variances with this algorithm are smaller in most cases. Meanwhile, it should be point out that when $l_1=l_2=0$, Algorithm \ref{AlgorithmSSCE} may give an inaccurate estimate, i.e., the ratio may be larger than $10$. This phenomenon also exists in estimating the normwise condition number of the LLS problem \cite{Babo}. Although the expression of $\kappa_{LSE} (A,B,b,d)$ is more complicated than that of the normwise condition number of the LLS problem and the circumstances on these two problems are different, we believe that the underlying reason should be the same; the reader can refer to \cite{Babo} for a detailed explanation.

\begin{table}[htp]
%\begin{table}[htp]
  \centering
  \caption{The efficiency of statistical condition estimates with $\kappa(A)=n^{l_1}$ and $\kappa(B)=s^{l_2}$}\label{Table1}
 \begin{tabular}{lc|c|cccccc}
\hline
 &  & $s^{l_2}$& \multicolumn{2}{c}{$l_2=0$} & \multicolumn{2}{c}{$l_2=3$} & \multicolumn{2}{c}{$l_2=5$}  \\
\hline
   & $n^{l_1}$ & &  mean & variance  & mean & variance  & mean & variance  \\
\hline
 $\|r\|_2=10^{-4}$   & $l_1=0$      &  $r_{ssce}$ &1.0296e+001 & 8.4423e-019 & 1.4428e+000 & 2.4025e-001 & 1.1939e+000 & 2.4425e-001 \\

                               &    & $r_{pce}$  &1.0000e+000 & 1.1164e-011 & 1.0002e+000 & 1.5433e-007 & 1.0000e+000 & 2.6382e-008 \\

                    & $l_1=3$   &  $r_{ssce}$ &1.1779e+000 & 2.5183e-001 & 1.4516e+000 & 2.7551e-001 & 1.2142e+000 & 2.5935e-001  \\

                               &    & $r_{pce}$ &1.0000e+000 & 2.8406e-012 & 1.0001e+000 & 1.3532e-007 & 1.0000e+000 & 3.7685e-008 \\
                 & $l_1=5$ &  $r_{ssce}$ &1.1038e+000 & 2.6402e-001 & 1.0224e+000 & 2.6509e-001 & 1.0392e+000 & 2.6771e-001 \\

                               &    & $r_{pce}$ &1.0000e+000 & 7.8330e-012 & 1.0000e+000 & 1.9491e-011 & 1.0000e+000 & 1.8920e-011 \\
\hline
$\|r\|_2=10^{0}$ & $l_1=0$ &  $r_{ssce}$ &1.0296e+001 & 7.4623e-011 & 1.4464e+000 & 2.8249e-001 & 1.1772e+000& 2.6635e-001 \\

                               &    & $r_{pce}$ &1.0000e+000 & 1.7148e-012 & 1.0002e+000 & 1.5495e-007 & 1.0000e+000 & 1.3253e-008 \\

                 & $l_1=3$ &  $r_{ssce}$ &1.1175e+000 & 2.6104e-001&  1.0350e+000 & 2.7807e-001 & 1.0804e+000 & 2.5401e-001  \\

                               &    & $r_{pce}$ &1.0000e+000 & 1.1119e-011 & 1.0000e+000 & 1.8794e-011 & 1.0000e+000 & 1.8791e-011\\
                 & $l_1=5$ &  $r_{ssce}$ &1.0836e+000 &  3.0292e-001&  1.0585e+000 & 2.9561e-001 & 1.0428e+000 & 2.5327e-001 \\

                               &    & $r_{pce}$ &1.0000e+000 & 1.9198e-011 & 1.0000e+000 & 1.8810e-011 & 1.0000e+000 & 1.9298e-011\\
\hline
$\|r\|_2=10^{4}$ & $l_1=0$ &  $r_{ssce}$ &9.4308e+000 & 5.8316e-003 & 1.4344e+000 & 2.4273e-001 & 1.2100e+000 & 2.6137e-001 \\

                               &    & $r_{pce}$ &1.0000e+000 & 6.0233e-014 & 1.0002e+000 & 2.0059e-007 & 1.0000e+000 & 2.5855e-008\\

                 & $l_1=3$ &  $r_{ssce}$ &1.0732e+000 & 2.8682e-001 & 1.0666e+000 & 2.8126e-001&  1.0120e+000 & 2.5842e-001 \\

                               &    & $r_{pce}$ &1.0000e+000 & 1.8471e-011 & 1.0000e+000 & 1.8770e-011 & 1.0000e+000 & 1.9124e-011\\
                 & $l_1=5$ &  $r_{ssce}$ &9.9248e-001 & 2.5030e-001 & 1.0454e+000 & 2.9353e-001&  1.0158e+000 & 2.9255e-001 \\

                               &    & $r_{pce}$ &1.0000e+000 & 1.9461e-011 & 1.0000e+000 & 1.9437e-011 & 1.0000e+000 & 1.9347e-011\\
\hline
\end{tabular}
%\end{table}
\end{table}

}\end{example}

\begin{example}{
Let $A$ and $B$ be gaussian random Toeplitz matrices of order $n=100$. This means that the entries of these two matrices are generated from standard normal distribution. Analogous to Example \ref{example1}, we also let the solution $x$ be $x=(1, 2^2,\cdots, n^2)^T$ and the residual vector $r$ be a random vector of specified norm. However, unlike Example \ref{example1}, it seems impossible to restrict a specific condition number to gaussian random Toeplitz matrices.%, but a good estimate of  its condition number  may be available \cite{Pan12}.

In our numerical experiment, for each $r$, we test $200$ pairs of random Toeplitz matrices $A$ and $B$. The numerical results on the ratio between $\kappa_{LSE} (A,B,b,d) $ and $\kappa_{LSE}^S (A,B,b,d)$ defined by
$$ratio= \frac{\kappa_{LSE} (A,B,b,d)}{\kappa_{LSE}^S (A,B,b,d)}$$
are presented in Figure \ref{Fig1}, which confirms the theoretical analysis in Remark \ref{rmk3.2}.
\begin{figure}[htbp]
  \centering
  \includegraphics[width=1\textwidth,height=0.85\textwidth]{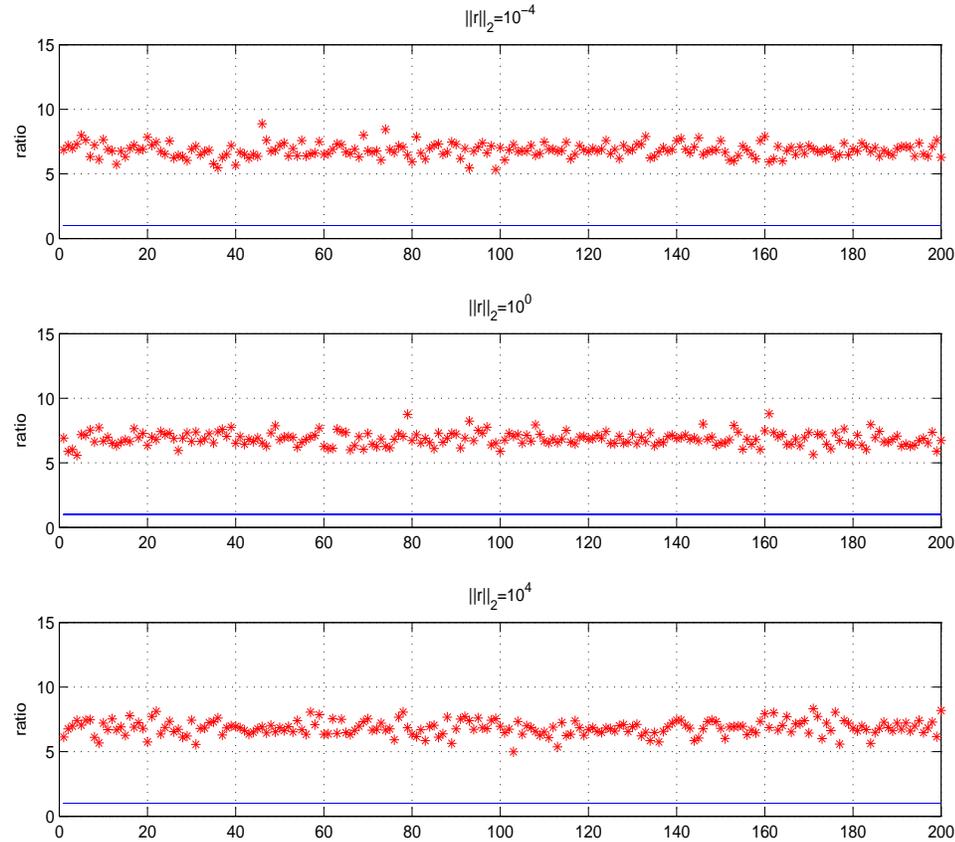}\\
  \caption{Comparison of $\kappa_{LSE} (A,B,b,d) $ and $\kappa_{LSE}^S (A,B,b,d)$ }\label{Fig1}
\end{figure}
}\end{example}
From Figure \ref{Fig1}, we also find that there are some points near $10$, which means that the unstructured condition number can be $10$ times larger than the structured one. Thus, it may lead to an overestimate when using the unstructured condition number to give error bounds in a structured LSE problem. Moreover, we also note that, for different $\|r\|_2$s, the $ratio$s seem to follow the same trend gathering in the interval $[5, 10]$. Whereas, from numerical experiments, we verify that the $ratio$ tends to be larger as $n$ increases. In the numerical experiments, we set  $\|r\|_2=1$ and $n=20*i-10,\; i=1:11$, and, for every $n$, we test $50$ LSE problems with random Toeplitz coefficient matrices  $A$ and $B$. The numerical results are presented in Figure \ref{Fig2}, where the circle line denotes the mean value of $ratio$s and the solid line denotes the corresponding variances. The fact shown in the figure means that the structured condition number has more advantage compared with the unstructured one as the dimensions of coefficient matrices increase.

% measuring the $ratio$'s variability in some sense.
\begin{figure}[htbp]
  \centering
  \includegraphics[width=0.81\textwidth,height=0.47\textwidth]{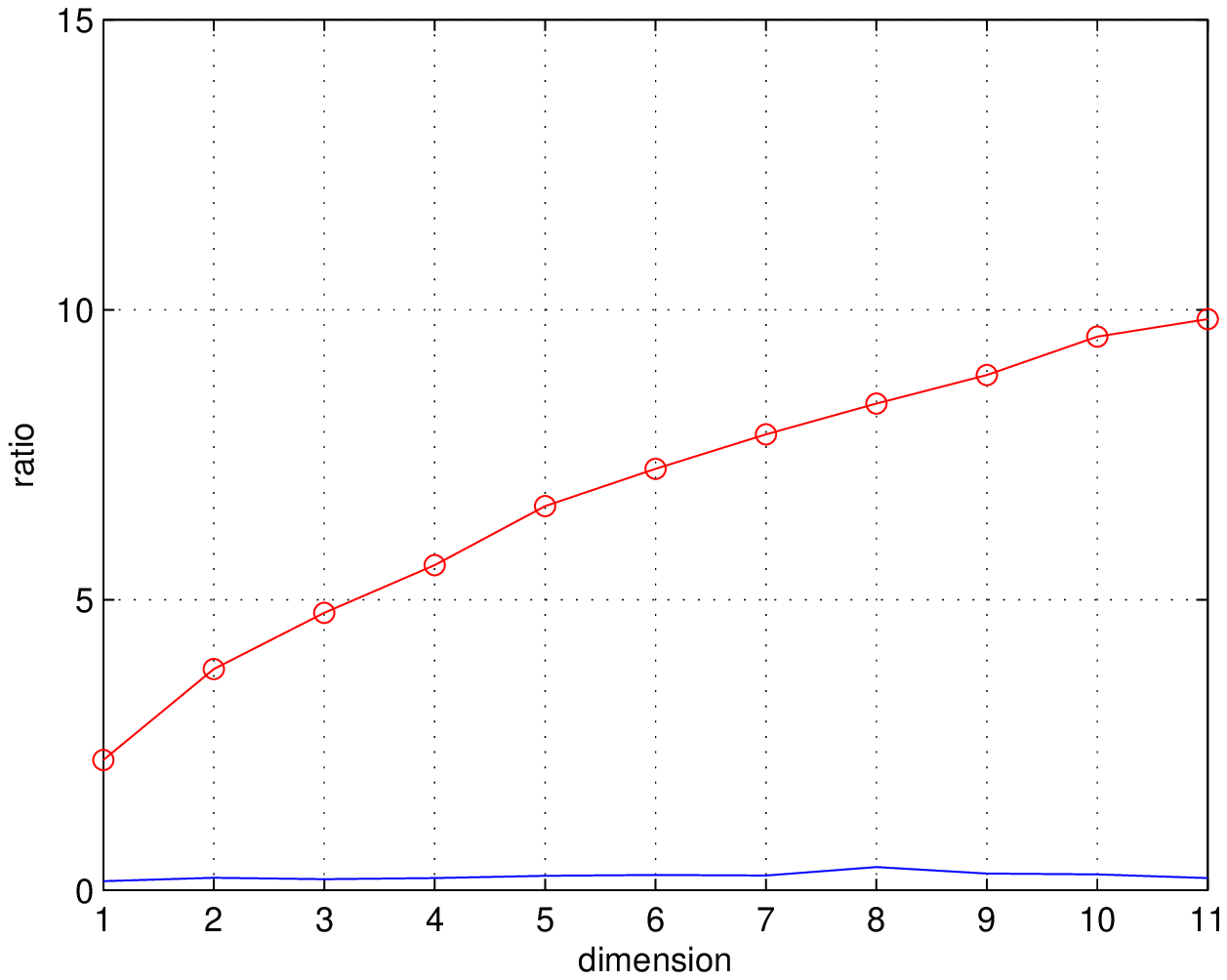}\\
  \caption{The influence of dimension}\label{Fig2}
\end{figure}

%\begin{remark}{\em
%From Figure \ref{Fig1}, we can see that all $ratio$s are located above the blue line, which shows that the all the $ratio$s are larger than $1$ and justifies our theoretical proof. In Figure \ref{Fig2}, we note that as the order of coefficient matrices increases the $ratio$ can be larger and larger, which means that the structured condition number will give tighter error bound when the dimensions of coefficient matrices are large compared with unstructured one. Due to the complicated expression of $\kappa_{LSE}^S (A,B,b,d)$, an appropriate explanation for the phenomenon is not easy to give, but in some sense this coincides with the results given in \cite{Pan12} that the condition number of gaussian random Toeplitz matrix becomes larger as the order increases.
%}\end{remark}

%\section{Concluding remarks}
%This paper presents the explicit expression and closed formulae of the partial condition number of the LSE problem. These results generalize the corresponding ones for the LLS problem given in \cite{Ari} and improve the one for the LSE problem in \cite{Cox}. The indefinite least squares problem with equality constraints is a generalization of the LSE problem \cite{Bob}, and finds some applications in the solution of total least squares problems and the area of optimization known as $H^\infty$ smoothing \cite{Bob,Liu101}. Recently, some authors considered its algorithms and mixed and componentwise condition numbers \cite{Liu101,Liu10,Mas1,Mas2,Li14}. It is also interesting to study its partial condition number. We will consider this topic in the future.

\section*{Acknowledgments}
{The authors would like to thank Prof. Michiel E. Hochstenbach for providing Matlab program of the probabilistic spectral norm estimator.}

\end{document}